# On a geometric black hole of a compact manifold

Alexander A. Ermolitski

*Cathedra of mathematics, MSHRC,
st. Nezavisimosti 62, Minsk, 220050, Belarus
E-mail: ermolitski@mail.by*

___

**Abstract:** Using a smooth triangulation and a Riemannian metric on a compact, connected, closed manifold $M^n$ of dimension $n$ we claim that every such $M^n$ can be represented as a union of a $n$–dimensional cell $C^n$ and a connected union $K^{n-1}$ (dim $K^{n-1} \leq n - 1$) of some finite number of subsimplexes of the triangulation. A sufficiently small closed neighborhood of $K^{n-1}$ is called a *geometric black hole*. Any sooth tensor field $K$ (or other structures) can be deformed into a continuous and sectionally smooth tensor field $\overline{K}$ where $\overline{K}$ has a very simple construction out of the black hole.

___

## 1. On extension of coordinate neighborhood

**1°.** Let $M^n$ be a connected, compact, closed and smooth manifold of dimension $n$ and $C^n$ be a cell (coordinate neighborhood) on $M^n$. We can fix some Riemannian metric $g$ on the manifold $M^n$ which defines the length of arc of a piecewise smooth curve and the continuous function $\rho(x;y)$ of the distance between two points $x, y \in M^n$. The topology defined by the function of distance (metric) $\rho$ is the same as the topology of the manifold $M^n$, [3].

A standard simplex $\Delta^n$ of dimension $n$ is the set of points $x=(x_1, x_2, ..., x_n) \in \mathbf{R}^n$ defined by conditions
$$0 \leq x_i \leq 1, \ i=\overline{1, \ n}, \ x_1+x_2+...+x_n \leq 1.$$

We consider the interval of a straight line connected the center of some face of $\Delta^n$ and the vertex which is opposite to this face. It is clear that the center of $\Delta^n$ belongs to the interval. We can decompose $\Delta^n$ as a set of intervals which are parallel to that mentioned above. If the center of $\Delta^n$ is connected by intervals with points of some face of $\Delta^n$ then a subsimplex of $\Delta^n$ is obtained. All the faces of $\Delta^n$ considered, $\Delta^n$ is seen as a set of all such subsimplexes. Let $U(\Delta^n)$ be some open neighborhood of $\Delta^n$ in $\mathbf{R}^n$. A diffeomorphism $\varphi : U(\Delta^n) \to M^n (\delta^n = \varphi(\Delta^n))$ is called a singular $n$–simplex on the manifold $M^n$. *Faces, edges, the center, vertexes* of the simplex $\delta^n$ are defined as the images of those of $\Delta^n$ with respect to $\varphi$.

The manifold $M^n$ is triangulable, **[5]**. It means that for any $l$, $0 \leq l \leq n$ such a finite set $\Phi^l$ of diffeomorphisms $\varphi : \Delta^l \to M^n$ is defined that

a) $M^n$ is a disjunct union of images $\varphi(Int\Delta^l)$, $\varphi \in \Phi^l$;

b) if $\varphi \in \Phi^l$ then $\varphi \circ \varepsilon_i \in \Phi^{l-1}$ for every $i$ where $\varepsilon_i : \Delta^{k-1} \to \Delta^k$ is the linear mapping transferring the vertexes $v_0,...,v_{k-1}$ of the simplex $\Delta^{k-1}$ in the vertexes $v_0,...,\hat{v}_i,...v_k$ of the simplex $\Delta^k$.

**2°.** Let $\delta_0^n$ be some simplex of the fixed triangulation of the manifold $M^n$. We paint the inner part $Int\delta_0^n$ of the simplex $\delta_0^n$ in white and the boundary $\partial\delta_0^n$ of $\delta_0^n$ in black. There exist coordinates on $Int\delta_0^n$ given by diffeomorphism $\varphi_0$. A subsimplex $\delta_{01}^{n-1} \subset \delta_0^n$ is defined by a black face $\delta_{01}^{n-1} \subset \delta_0^n$ and the center $c_0$ of $\delta_0^n$. We connect $c_0$ with the center $d_0$ of the face $\delta_{01}^{n-1}$ and decompose the subsimplex $\delta_{01}^n$ as a set of intervals which are parallel to the interval $c_0d_0$. The face $\delta_{01}^{n-1}$ is a face of some simplex $\delta_1^n$ that has not been painted. We draw an interval between $d_0$ and the vertex $v_1$ of the subsimplex $\delta_1^n$ which is opposite to the face $\delta_{01}^{n-1}$ then we decompose $\delta_1^n$ as a set of intervals which are parallel to the interval $d_0v_1$. The set $\delta_{01}^n \cup \delta_1^n$ is a union of such broken lines every one from which consists of two intervals where the endpoint of the first interval coincides with the beginning of the second interval (in the face $\delta_{01}^{n-1}$) the first interval belongs to $\delta_{01}^n$ and the second interval belongs to $\delta_1^n$. We construct a homeomorphism (extension) $\varphi_{01}^1 : Int\delta_{01}^n \to Int(\delta_{01}^n \cup \delta_1^n)$. Let us consider a point $x \in Int\delta_{01}^n$ and let $x$ belong to a broken line consisting of two intervals the first interval is of a length of $s_1$ and the second interval is of a length of $s_2$ and let $x$ be at a distance of $s$ from the beginning of the first interval. Then we suppose that $\varphi_{01}^1(x)$ belongs to the same broken line at a distance of $\dfrac{s_1+s_2}{s_1} \cdot s$ from the beginning of the first interval. It is clear that $\varphi_{01}^1$ is a homeomorphism giving coordinates on $Int(\delta_{01}^n \cup \delta_1^n)$. We paint points of $Int(\delta_{01}^n \cup \delta_1^n)$ white. Assuming the coordinates of points of white initial faces of subsimplex $\delta_{01}^n$ to be fixed we obtain correctly introduced coordinates on $Int(\delta_0^n \cup \delta_1^n)$. The set $\sigma_1 = \delta_0^n \cup \delta_1^n$ is called a *canonical polyhedron*. We paint faces of the boundary $\partial\sigma_1$ black.

We describe the contents of the successive step of the algorithm of extension of coordinate neighborhood. Let us have a canonical polyhedron $\sigma_{k-1}$ with white inner points (they have introduced *white coordinates*) and the black boundary $\partial \sigma_{k-1}$. We look for such an *n*–simplex in $\sigma_{k-1}$, let it be $\delta_0^n$ that has such a black face, let it be $\delta_{01}^{n-1}$ that is simultaneously a face of some *n*–simplex, let it be $\delta_1^n$, inner points of which are not painted. Then we apply the procedure described above to the pair $\delta_0^n$, $\delta_1^n$. As a result we have a polyhedron $\sigma_k$ with one simplex more than $\sigma_{k-1}$ has. Points of $Int \sigma_k$ are painted in white and the boundary $\partial \sigma_k$ is painted in black. The process is finished in the case when all the black faces of the last polyhedron border on the set of white points (the cell) from two sides.

After that all the points of the manifold $M^n$ are painted in black or white, otherwise we would have that $M^n = M_0^n \cup M_1^n$ (the points of $M_0^n$ would be painted and those of $M_1^n$ would be not) with $M_0^n$ and $M_1^n$ being unconnected, which would contradict of connectivity of $M^n$.

Thus, we have proved the following

**Theorem 1.** *Let $M^n$ be a connected, compact, closed, smooth manifold of dimension n. Then $M^n = C^n \cup K^{n-1}$, $C^n \cap K^{n-1} = \varnothing$, where $C^n$ is an n–dimensional cell and $K^{n-1}$ is a union of some finite number of (n–1)–simplexes of the triangulation.*

3°. We consider the initial simplex $\delta_0^n$ of the triangulation and its center $c_0$. Drawing intervals between the point $c_0$ and points of all the faces of $\delta_0^n$ we obtain a decomposition of $\delta_0^n$ as a set of the intervals. In **2°** the homeomorphism $\psi: Int \delta_0^n \to C^n$ was constructed and $\psi$ evidently maps every interval above on a piecewise smooth broken line $\gamma$ in $C^n$. We denote $\widetilde{M}^n = M^n \setminus \{c_0\}$. $\widetilde{M}^n$ is a connected and simply connected manifold if $M^n$ is that. Let $I=[0;1]$, we define a homotopy $F: \widetilde{M}^n \times I \to \widetilde{M}^n : (x; t) \mapsto y = F(x;t)$ in the following way

a) $F(z; t) = z$ for every point $z \in K^{n-1}$;

b) if a point $x$ belongs to the broken line $\gamma$ in $C^n$ and the distance between $x$ and its limit point $z \in K^{n-1}$ is $s(x)$ then $y = F(x; t)$ is on the same broken line $\gamma$ at a distance of $(1-t)s(x)$ from the point $z$.

It is clear that $F(x;0) = x$, $F(x;1) = z$ and we have obtained the following

**Theorem 2.** *The spaces $\widetilde{M}^n$ and $K^{n-1}$ are homotopy–equivalent, in particular, the groups of singular homologies $H_k(\widetilde{M}^n)$ and $H_k(K^{n-1})$ are isomorphic for every k.*

**Corollary 2.1.** *The space $K^{n-1}$ is connected and if $M^n$ is simply connected then $K^{n-1}$ is simply connected too.*

**Remark.** *The white coordinates are extended from the simplex $\delta_0^n$ in the simplex $\delta_1^n$ through the face $\delta_{01}^{n-1}$ hence $Int \delta_{01}^{n-1}$ has also the white coordinates.*

*On the other hand there exist two linear structures (intervals, the center etc) on $\delta_{01}^n$ induced from $\delta_0^n$ and $\delta_1^n$ respectively. Further, we set that the linear structure of $\delta_{01}^{n-1}$ is the structure induced from $\delta_0^n$.*

## 2. Deformation of a tensor field towards a geometric black hole of a compact manifold

**1°.** Let $L(M^n)$ and $L(C^n)$ be the principal fibre bundles of linear frames of the manifolds $M^n$ and $C^n$. The diffeomorphism $\varphi_0$ (where $\delta_0^n = \varphi_0(\Delta^n)$) defines the coordinates $(x_1, ..., x_n)$ in some neighborhood of the simplex $\delta_0^n$ and the corresponding vector fields $\frac{\partial}{\partial x_1},...,\frac{\partial}{\partial x_n}$ on this neighborhood (a local cross–section of $L(M^n)$). Similarly, the diffeomorphism $\varphi_1$ (where $\delta_1^n = \varphi_1(\Delta^n)$) defines the coordinates $(y_1, ..., y_n)$ in some neighborhood of the simplex $\delta_1^n$ and the vector fields $\frac{\partial}{\partial y_1},...,\frac{\partial}{\partial y_n}$ on this neighborhood. We have assumed that the white face $\delta_{01}^{n-1}$ has the equation: $y_1=0$ (it can always be obtained by corresponding linear change of variables in $\mathbf{R}^n \supset \Delta^n$). The vector fields $X_i = \frac{\partial}{\partial x_i}, \frac{\partial}{\partial y_j}$, $i,j = \overline{1,n}$, are defined on the face $\delta_{01}^{n-1}$ therefore for any point $x \in \delta_{01}^{n-1}$ we have $X_i = \sum_{j=1}^{n} f_{ij}(x)\frac{\partial}{\partial y_i}$ where the functions $f_{ij}(x)$ are smooth. We decompose $\delta_1^n$ as a union of the intervals having the following equations: $y_1=t$, $y_2=c_2$, $y_3=c_3,..., y_n=c_n$, where $0, c_2,..., c_n$ are the coordinates of the beginning $y_0$ of the corresponding interval. For any point $y \in \delta_1^n$ we assume $f_{ij}(y) = f_{ij}(y_0)$ where $y_0 \in \delta_{01}^{n-1}$ is the beginning of the interval where the point $y$ is situated. The vector fields $X_i$, $i = \overline{1,n}$, are defined on $\delta_1^n$ by the formula $X_i = \sum_{j=1}^{n} f_{ij}(y)\frac{\partial}{\partial y_i}$. It is obvious that the constructed vector fields $X_i$, $i = \overline{1,n}$, are continuous on $\delta_0^n \cup \delta_1^n$ and smooth in any point $x \in \delta_0^n \cup \delta_1^n$, $x \notin \delta_{01}^{n-1}$.

For the process of the extension of a coordinate neighborhood (**1, 2°**) we can consider the process of the extension of the vector fields $X_1,..., X_n$. If these fields are defined on a polyhedron $\sigma_{k-1}$ and in order to get a polyhedron $\sigma_k$ we use simplexes $\delta_0^n, \delta_1^n$ then we apply the procedure described above to obtain the vector fields on $\sigma_k$. As a result we obtain correctly defined vector fields $X_1,..., X_n$ on $C^n$ i.e. a cross–section of $L(C^n)$.

So, we come to the following

**Proposition 3.** *Let $M^n = C^n \cup K^{n-1}$ be the decomposition from the theorem 1. Then there exists a continuous cross–section of $L(C^n)$: $x \to (X_1,\ldots, X_n)_x$, $x \in C^n$. If a point $x \in C^n$ does not belong to the subsimplexes of the triangulation then the cross–section above is smooth at the point x.*

We consider a tensor of type $(r, s)$ on $\mathbf{R}^n$:

$$K^0 = \sum k^{\lambda_1\ldots\lambda_r}_{\mu_1\ldots\mu_s}(0) e_{\lambda_1} \otimes \ldots \otimes e_{\lambda_r} \otimes e^{\mu_1} \otimes \ldots \otimes e^{\mu_s},$$

where $e_1,\ldots,e_n$ is the standard basic of $\mathbf{R}^n$ and $e^1,\ldots,e^n$ is the dual basis of $\mathbf{R}^{n*}$.

A tensor field of type $(r, s)$ is defined on $C^n$:

$$K^0 = \sum k^{\lambda_1\ldots\lambda_r}_{\mu_1\ldots\mu_s}(0) X_{\lambda_1} \otimes \ldots \otimes X_{\lambda_r} \otimes X^{\mu_1} \otimes \ldots \otimes X^{\mu_s} \tag{1}$$

Since the functions $k^{\lambda_1\ldots\lambda_r}_{\mu_1\ldots\mu_s}$ are constant on $C^n$ we obtain that the tensor field $K^0$ is O–deformable on $C^n$ i.e. some G–structure on $C^n$ is defined by $K^0$ (see [1], [4]). If the cross–section $(X_1,\ldots, X_n)_x$ is smooth at a point $x \in C^n$ then the tensor field $K^0$ is also smooth at the point.

**2°.** For any point $z \in K^{n-1}$ we can consider the closed geodesic ball $\overline{B}(z,\varepsilon)$ of a small radius $\varepsilon > 0$. Let $Tb(K^{n-1}, \varepsilon) = \bigcup_{z \in K^{n-1}} \overline{B}(z,\varepsilon) = BH(\varepsilon)$.

**Definition 1.** *We call the set $BH(\varepsilon)$ a geometric black hole of radius $\varepsilon>0$ of the manifold $M^n$ if $M^n \setminus BH(\varepsilon)$ is cell (it is true for some small $\varepsilon$). We paint the points of $BH(\varepsilon)$ in black.*

Any piecewise smooth broken line $\gamma$ considered in **1, 3°** can be represented as $\gamma = \gamma_0 \cup \gamma_1$ where $\gamma_1 = \gamma \cap BH(\varepsilon)$, $\gamma_0 = \gamma \setminus \gamma_1$. The points of $\gamma_0$ are painted in white and the points of $\gamma_1$ are painted in black. Let the segment $\gamma_0$ have a length $s_0$ and the segment $\gamma_1$ have a length $s_1$ then $(s_0+s_1)$ is a length of the broken line $\gamma$ from $c_0$ to $z \in K^{n-1}$.

Let $K(x)$, $x \in M^n$, be a tensor field of type $(r, s)$ and $K^0 = K(c_0)$ where $c_0$ is the center of the initial simplex $\delta^n_0$ of the triangulation of $M^n$. Deformations of structures was considered in [2]. So, we shall construct a deformation $\overline{K}(x)$ of the tensor field $K(x)$ on the manifold $M^n$.

1) If a point $z \in K^{n-1}$ then $\overline{K}(z) = K(z)$.

2) If a point $x \in M^n \setminus BH(\varepsilon)$ then $\overline{K} = K^0 = K(c_0)$ where $K^0$ is defined by the formula (1).

3) *We assume that* $K(x) = \sum k^{\lambda_1\ldots\lambda_r}_{\mu_1\ldots\mu_s}(x) X_1 \otimes \ldots \otimes X_{\lambda_r} \otimes X^{\mu_1} \otimes \ldots \otimes X^{\mu_s}$,
$x \in C^n$, where $X_1,\ldots, X_r$ are the vector fields from the proposition 3, a point $x$ belongs a broken line $\gamma$ and $s(x)$ is the distance from $x$ to $c_0$ along the broken line $\gamma$. For any point $y \in \gamma_1$ we define the tensor field

$$\overline{K}(y) = \sum \overline{k}^{\lambda_1\ldots\lambda_r}_{\mu_1\ldots\mu_s}(y) X_1 \otimes \ldots \otimes X_{\lambda_r} \otimes X^{\mu_1} \otimes \ldots \otimes X^{\mu_s}$$

in the following way: $\bar{k}^{\lambda_1\ldots\lambda_r}_{\mu_1\ldots\mu_s}(y) = k^{\lambda_1\ldots\lambda_r}_{\mu_1\ldots\mu_s}(x)$ where $s(x) = \frac{s(y) - s_0}{s_1}(s_0 + s_1)$, $s(y)$ is the distance from $y$ to $c_0$ along the broken line $\gamma$.

It is easy to see that the constructed tensor field $\bar{K}$ is continuous and sectionally smooth, $\bar{K}$ is not smooth on the boundary of $BH(\varepsilon)$ and in the points of $C^n$ where the cross–section $(X_1,\ldots, X_n)_x$ is not smooth.

**Remark** *We can also consider similar deformations for other geometric structures towards* $BH(\varepsilon)$.

# References

[1] A.A. Ermolitski, *Riemannian manifolds with geometric structures,* BSPU, Minsk, 1998 (in Russian), (English version: Internet, Google, Ermolitski).
[2] A.A. Ermolitski, *Deformations of structures, embedding of a Riemannian manifold in a Kaehlerian one and geometric antigravitation,* Banach Center Publications, Warszawa, 2007, V. 76, 505–514.
[3] D. Gromoll, W. Klingenberg, W. Meyer, *Riemannsche Geometrie im Grossen,* Springer, Berlin, 1968.
[4] S. Kobayashi, *Transformation groups in differential geometry,* Springer, Berlin–New York, 1972.
[5] J.R. Munkres, *Elementary differential topology,* Princeton University Press, Princeton, 1966.